\documentclass{ifacconf}

\usepackage{graphicx}      
\usepackage{natbib}        
\usepackage{amsfonts}
\usepackage{amsmath}
\usepackage{xcolor}

\newcommand{\diag}{\mathop{\mathrm{diag}}}

\allowdisplaybreaks

\begin{document}
\begin{frontmatter}

\title{Stability Properties of Multi-Order Fractional Differential Systems in 3D} 


\author[First]{Kai Diethelm} 
\author[First]{Safoura Hashemishahraki} 
\author[Second]{Ha Duc Thai}
\author[Second]{Hoang The Tuan}

\address[First]{Faculty of Applied Sciences and Humanities, 
Technical University of Applied Sciences Würzburg-Schweinfurt, Ignaz-Schön-Str.~11,
97421~Schweinfurt, Germany\\
(e-mail: \{kai.diethelm, safoura.hashemishahraki\}@thws.de)}
\address[Second]{Institute of Mathematics, Vietnam Academy of Science and Technology, 
18 Hoang Quoc Viet, 10307 Ha Noi, Vietnam\\ (e-mail: \{hdthai, httuan\}@math.ac.vn)}

\begin{abstract}                
This paper is devoted to studying three-dimensional non-commensurate fractional order differential equation systems with Caputo derivatives. Necessary and sufficient conditions are for the asymptotic stability of such systems are obtained. 
\end{abstract}


\end{frontmatter}

\section{Introduction}
Ordinary differential equation systems of fractional order with Caputo operators have been a topic of great interest for the last decades
because of their considerable importance in many applications in various fields of applied sciences. 
One relevant aspect when discussing such systems is their stability analysis. 
While the investigation of commensurate systems (i.e.\ systems where all equations contain differential operators of the same order)
is straightforward and well investigated \citep{Ma1996,Ma1998}, much less is known in the general incommensurate case where all equations are allowed to have arbitrary orders, independent of the orders of the other equations of the system. 
In this work we shall concentrate on the latter case and establish a few new results that provide additional insight.

The fundamental result about $d$-dimensional \emph{linear} systems of this type with some constant coefficient matrix $A$ 
is that the system is asymptotically stable if and only if the 
zeros of its fractional characteristic function $\det(\diag(s^{\alpha_1}, \ldots, s ^{\alpha_d}) - A)$ are in the open left half of the complex plane, see \cite{DengLiLu}. While this result is very valuable theoretically, it is only of rather limited practical use because finding roots of a fractional characteristic polynomial of a non-commensurate fractional-order system is a complicated task that has so far only been solved in some special cases, essentially in two dimensions; cf.\ \cite{BrandiburKaslik} and \cite{DiethelmThaiTuan} where, in particular, the results have also been extended to some classes of nonlinear equations.

Our aim here is to extend the work of those papers to three dimensions. The general $d$-dimensional setting will be discussed elsewhere later. We thus study here the fractional-order system with Caputo fractional derivatives 
\begin{eqnarray*}
	{D}^{\alpha} x(t) & = & A x(t) + f(t, x(t)), \,\,\, t > 0, \\
	                   x(0) & = & x^{0} \in \mathbb{R}^{3},
\end{eqnarray*}
where $\alpha = (\alpha_{1}, \alpha_{2}, \alpha_{3}) \in (0, 1]^3$ is a multi-index, $A \in \mathbb{R}^{3 \times 3}$ is a square matrix and $f : [0, \infty) \times \mathbb{R}^{3} \rightarrow \mathbb{R}^{3}$ is a vector valued continuous function. Our results complement those recently developed for general $d$ by \cite{DiethelmHashemishahrakiThaiTuan} where completely different methods are used and stability is shown under significantly different assumptions.

We will first recall some basic general properties of inhomogeneous linear systems of the type of equations that we are interested in. Next we will look at the concrete question for the asymptotic stability or instability of homogeneous linear systems, in particular providing criteria that can be effectively tested for a given system.
Afterwards, we show how these criteria can also be exploited to investigate inhomogeneous systems and systems with certain types of nonlinearities. Some examples will conclude our work.

\section{Some general properties of non-commensurate systems and their characteristic function}

Consider the three-dimensional incommensurate system
\begin{equation}
\label{eq:linearsystem}
\left.
  \begin{array}{ll}
    D^{\alpha_1} x_{1}(t) & \! = \! a_{11} x_{1}(t) + a_{12} x_{2} (t) + a_{13} x_{3}(t) + f_{1}(t),  \\
    D^{\alpha_2} x_{2}(t) & \! = \! a_{21} x_{1}(t) + a_{22} x_{2} (t) + a_{23} x_{3}(t) + f_{2}(t),  \\
    D^{\alpha_3} x_{3}(t) & \! = \! a_{31} x_{1}(t) + a_{32} x_{2} (t) + a_{33} x_{3}(t) + f_{3}(t),
  \end{array}
  \!
\right\}
\end{equation}
where $x_{1}(0) = x_{1}^{0}$, $x_{2}(0) = x_{2}^{0}$, $x_{3}(0) = x_{3}^{0}$
and  $\alpha_{1},\alpha_{2},\alpha_{3} \in (0,1]$.
We use the obvious notation
\begin{equation}
A = \begin{pmatrix}
        a_{11} & a_{12} & a_{13} \\
        a_{21} & a_{22} & a_{23} \\
        a_{31} & a_{32} & a_{33} 
      \end{pmatrix},
      \quad
       f (t) = \begin{pmatrix}
                          f_1(t) \\
                          f_2(t) \\
                          f_3(t) 
                        \end{pmatrix}.
\end{equation}
By applying the Laplace transform to \eqref{eq:linearsystem}, we obtain
\begin{equation*}
	A \cdot \begin{pmatrix}
                 X_1 (s) \\
                 X_2 (s) \\
                 X_3 (s) 
               \end{pmatrix}
      = \begin{pmatrix}
                           s^{\alpha_{1}-1} x_{1}^0 + F_1(s) \\
                           s^{\alpha_{2}-1} x_{2}^0 + F_2(s) \\
                           s^{\alpha_{3}-1} x_{3}^0 + F_3(s)
          \end{pmatrix}
\end{equation*}
where $X_k$ and $F_k$ denote the Laplace transforms of $x_k$ and $f_k$, respectively. 
Then we can find the solution $(X_1(s), X_2(s), X_3(s))^{\mathrm T}$ of this algebraic linear system with standard linear algebra methods. Finally an inverse Laplace transform yields the desired solution $(x_1(t), x_2(t), x_3(t))^{\mathrm T}$ of the differential equation system~\eqref{eq:linearsystem}.

Classical Laplace transform arguments then show us the following result, see \cite{DengLiLu}:

\begin{thm}
\label{thm:deng}
The homogeneous system associated to \eqref{eq:linearsystem} is asymptotically stable if all zeros of its characteristic function $Q$ 
defined by
\[
	Q(s) = \det \left( \begin{pmatrix}
						s^{\alpha_1} & 0 & 0 \\
						0 & s^{\alpha_2} & 0 \\
						0 & 0 & s^{\alpha_3} 
						\end{pmatrix} - A \right)
\]
are in the open left half of the complex plane. Moreover, the homogeneous system is unstable if $Q$ has at least one zero in the open right half of $\mathbb C$.
\end{thm}

Therefore, our main question is:
Where are the zeros of the characteristic function $Q$?
More precisely, we are interested in finding conditions on the matrix $A$ which are sufficient to guarantee that all zeros of $Q$ are in the open left half of the complex plane.

An explicit calculation shows that the characteristic function $Q$ has the form
\begin{eqnarray}
Q(s) & = & s^{\alpha_1 + \alpha_2 + \alpha_3} - a_{11} s^{\alpha_2 + \alpha_3} - a_{22} s^{\alpha_1 + \alpha_3} - a_{33} s^{\alpha_1 + \alpha_2} 
\nonumber\\
 && {} + s^{\alpha_{1}} (a_{22} a_{33}-a_{23} a_{32}) + s^{\alpha_{2}} (a_{11} a_{33} - a_{13} a_{31})
 \nonumber\\
 && {} + s^{\alpha_{3}} (a_{11} a_{22} - a_{12} a_{21}) - \det A.
\label{eq:q}
\end{eqnarray}
Our considerations below will be based on this relation.

\section{Main Results}
\label{sec:mainresults}

Following the standard approach, we will first discuss the homogeneous linear system associated
to~\eqref{eq:linearsystem}, i.e.\ the case where $f_k(t) = 0$ for all $k$.
We then extend these results to the inhomogeneous case (see Subsection \ref{subs:inhom-lin})
and consider nonlinear problems (cf.\ Subsection \ref{subs:nonlin}).

\subsection{Homogeneous linear systems}
\label{subs:lin-hom}

Our first new theorem is a general negative result.

\begin{thm}
	\label{thm:unstable1}
	\begin{enumerate}
	\item[(a)] If $\det A > 0$ then the homogeneous system associated to \eqref{eq:linearsystem} is unstable.
	\item[(b)] If $\det A < 0$, $\alpha_1 + \alpha_2 + \alpha_3 \ge 2$ and $Q$ is of the form 
		$Q(s) = s^{\alpha_1 + \alpha_2 + \alpha_3} - \sum_{k=1}^n c_k s^{\psi_k} - \det A$
		with some $n \in \mathbb N$, $\max_k c_k > \min_k c_k \ge 0$ and $\psi_k \in (0,2)$ 
		then the homogeneous system associated to \eqref{eq:linearsystem} is unstable.
	\end{enumerate}
\end{thm}

\begin{pf}
	From eq.~\eqref{eq:q}, we can see 
	that $Q$ is continuous on $\mathbb R$, that $\lim_{s \to \infty} Q(s) = \infty$,
	and that $Q(0) = -\det A$.
	
	Hence, if $\det A > 0$ then $Q(0) < 0$ and thus $Q$ must have a real zero in the interval $(0,\infty)$.
	Therefore, the statement of Theorem~\ref{thm:deng} asserts the instability in case (a). 
	
	We will give the proof of (b) at the end of Subsection \ref{subs:lin-hom}.~\qed
\end{pf}

\begin{rem}
	Among all the results that we shall prove in this paper, Theorem \ref{thm:unstable1}
	is special because it is the only statement that admits an immediate generalization from
	the three-dimensional setting to the general $d$-dimensional case with arbitrary $d \in \mathbb N$.
	To be precise, when following the same steps as above, one can see that the absolute term in the
	characteristic function $Q(s)$ has the value $(-1)^d \det A$. Therefore, arguing in exactly the same 
	manner as in the proof of Theorem \ref{thm:unstable1}, we can see that the homogeneous 
	system is unstable if $(-1)^d \det A < 0$ and not asymptotically stable if $\det A = 0$.
\end{rem}

Already in the two-dimensional case discussed in \cite{DiethelmThaiTuan} (where the coefficient matrix $A$
contains only four entries), it was necessary to discuss a relatively large number of cases in each of which the
coefficients were assumed to be in certain ranges. In the three-dimensional case that we discuss here, we have
nine entries in the matrix $A$, and therefore the number of cases is much bigger. Thus, space limitations
prevent us from discussing the matter exhaustively here. We therefore exclusively concentrate for the 
remainder of the paper on a limited range
of possible matrices, namely on those for which the characteristic function takes the simpler form
\begin{equation}
	\label{eq:char-func-simple}
	\left.
	\begin{array}{rl}
		Q(s) & = s^{\beta_4} - a s^{\beta_3} - b s^{\beta_2} - c s^{\beta_1} - d \\
		\text{ with some }
		0 & < \beta_1 \le \beta_2 \le \beta_3 < \beta_4
	\end{array}
	\right\}
\end{equation}
where, in particular, $\beta_4 = \alpha_1 + \alpha_2 + \alpha_3$
and $d = \det A$, i.e.\ we assume that three of the
six coefficients of the general form \eqref{eq:q} (not counting the leading coefficient, which always is $1$,
and the absolue term which is always $\det A$ and which, according to Theorem \ref{thm:unstable1},
must not be zero if asymptotic stability is desired) vanish. Our results below will therefore mainly
provide conditions under which the zeros of the function $Q$ of \eqref{eq:char-func-simple}
are located in the open left half of $\mathbb C$. Clearly, there are $\binom 6 3 = 20$ possible ways to select
the three coefficients that we assume to vanish from the six coefficients that are allowed to be zero in
principle. To be precise, our approach covers the following special situations where we always assume 
that the three mutually distinct values $j_1, j_2, j_3 \in \{ 1,2,3 \}$ are chosen such that 
$\alpha_{j_1} \le \alpha_{j_2} \le \alpha_{j_3}$:
\begin{enumerate}
\item $a_{11} = a_{22} = a_{33} = 0$; then $\beta_k = \alpha_{j_k}$ for $k = 1, 2, 3$;
\item $a_{11} = a_{22} = 0$ and $a_{22} a_{33} = a_{23} a_{32}$; then $\beta_1 = \min(\alpha_2, \alpha_3)$, $\beta_3 = \max (\alpha_{1} + \alpha_{2}, \alpha_{3})$, $\beta_2 = \alpha_{1} + 2 \alpha_{2} + \alpha_3 - \beta_1 - \beta_3$;
\item $a_{11} = a_{22} = 0$ and $a_{11} a_{33} = a_{13} a_{31}$; then $\beta_1 = \min(\alpha_1, \alpha_3)$, $\beta_3 = \max (\alpha_{1} + \alpha_{2}, \alpha_{3})$, $\beta_2 = 2 \alpha_{1} + \alpha_{2} + \alpha_3 - \beta_1 - \beta_3$;
\item $a_{11} = a_{22} = 0$ and $a_{11} a_{22} = a_{12} a_{21}$; then $\beta_1 = \min(\alpha_1, \alpha_2)$, $\beta_2 = \max(\alpha_1, \alpha_2)$, $\beta_3 = \alpha_1 + \alpha_2$;
\item $a_{11} = a_{33} = 0$ and $a_{22} a_{33} = a_{23} a_{32}$; then $\beta_1 = \min (\alpha_{2} , \alpha_{3})$, $\beta_{3} = \max(\alpha_{1} + \alpha_{3} , \alpha_{2})$, $\beta_2 = \alpha_1 + \alpha_{2} + 2 \alpha_{3} - \beta_{1} - \beta_{3}$;
\item $a_{11} = a_{33} = 0$ and $a_{11} a_{33} = a_{13} a_{31}$; then $\beta_1 = \min (\alpha_{1} , \alpha_{3})$, $\beta_2 = \max(\alpha_{1} , \alpha_{3})$, $\beta_3 = \alpha_{1} + \alpha_{3}$;
\item $a_{11} = a_{33} = 0$ and $a_{11} a_{22} = a_{12} a_{21}$; then $\beta_1 = \min (\alpha_{1} , \alpha_{2})$, $\beta_3 = \max(\alpha_{1} + \alpha_{3} , \alpha_{2})$, $\beta_2 = 2 \alpha_{1} +  \alpha_{2} + \alpha_{3} - \beta_{1} - \beta_{3}$;
\item $a_{22} = a_{33} = 0$ and $a_{22} a_{33} = a_{23} a_{32}$; then $\beta_1 = \min (\alpha_{2} , \alpha_{3})$, $\beta_2 = \max(\alpha_{2} , \alpha_{3})$, $\beta_3 = \alpha_{2} + \alpha_{3}$;
\item $a_{22} = a_{33} = 0$ and $a_{11} a_{33} = a_{13} a_{31}$; then $\beta_1 = \min (\alpha_{1} , \alpha_{3})$, $\beta_3 = \max(\alpha_1, \alpha_{2} + \alpha_{3})$, $\beta_2 = \alpha_{1} +  \alpha_{2} + 2 \alpha_{3} - \beta_{1} - \beta_{3}$;
\item $a_{22} = a_{33} = 0$ and $a_{11} a_{22} = a_{12} a_{21}$; then $\beta_1 = \min (\alpha_{1} , \alpha_{2})$, $\beta_3 = \max(\alpha_1, \alpha_{2} + \alpha_{3})$, $\beta_2 = \alpha_{1} + 2 \alpha_{2} + \alpha_{3} - \beta_{1} - \beta_{3}$;
\item $a_{11} = 0$, $a_{22} a_{33} = a_{23} a_{32}$ and $a_{11} a_{33} = a_{13} a_{31}$; then $\beta_{1} = \min (\alpha_{1} + \alpha_{2} , \alpha_{3})$, $\beta_3 = \max (\alpha_{1} + \alpha_{2} , \alpha_{1} + \alpha_{3})$, $\beta_2 = 2 \alpha_1 + \alpha_2 + 2 \alpha_3 - \beta_1 - \beta_3$;
\item $a_{11} = 0$, $a_{11} a_{22} = a_{12} a_{21}$ and $a_{11} a_{33} = a_{13} a_{31}$; then $\beta_1 = \alpha_{1}$, $\beta_2 = \min (\alpha_{1} + \alpha_{2} , \alpha_{1} + \alpha_{3})$, $\beta_{3} = \max (\alpha_{1} + \alpha_{2} , \alpha_{1} + \alpha_{3})$;
\item $a_{11} = 0$, $a_{11} a_{22} = a_{12} a_{21}$ and $a_{22} a_{33} = a_{23} a_{32}$; then $\beta_1 = \min (\alpha_{1} + \alpha_{3} , \alpha_{2})$, $\beta_3 = \max (\alpha_{1} + \alpha_{2} , \alpha_{1} + \alpha_{3})$, $\beta_2 = 2 \alpha_{1} + 2 \alpha_{2} + \alpha_{3} - \beta_{1} - \beta_{3}$;
\item $a_{22} = 0$, $a_{22} a_{33} = a_{23} a_{32}$ and $a_{11} a_{33} = a_{13} a_{31}$; then $\beta_1 = \min (\alpha_{1} + \alpha_{2} , \alpha_{3})$, $\beta_3 = \max (\alpha_{2} + \alpha_{3} , \alpha_{1} + \alpha_{2})$, $\beta_2 = \alpha_{1} + 2 \alpha_{2} + 2 \alpha_{3} - \beta_{1} - \beta_{3}$;
\item $a_{22} = 0$, $a_{11} a_{22} = a_{12} a_{21}$ and $a_{11} a_{33} = a_{13} a_{31}$; then $\beta_{1} = \min (\alpha_{2} + \alpha_{3} , \alpha_{1})$, $\beta_3 = \max (\alpha_{1} + \alpha_{2} , \alpha_{2} + \alpha_{3})$, $\beta_2 = 2 \alpha_{1} + 2 \alpha_{2} + \alpha_{3} - \beta_{1} - \beta_{3}$;
\item $a_{22} = 0$, $a_{11} a_{22} = a_{12} a_{21}$ and $a_{22} a_{33} = a_{23} a_{32}$; then $\beta_1 = \alpha_2$, $\beta_2 = \min (\alpha_{1} + \alpha_{2} , \alpha_{2} + \alpha_{3})$, $\beta_3 = \max (\alpha_{2} + \alpha_{3} , \alpha_{1} + \alpha_{2})$;
\item $a_{33} = 0$, $a_{22} a_{33} = a_{23} a_{32}$ and $a_{11} a_{33} = a_{13} a_{31}$; then $\beta_1 = \alpha_3$, $\beta_2 = \min (\alpha_{2} + \alpha_{3} , \alpha_{1} + \alpha_{3})$, $\beta_3 = \max (\alpha_{1} + \alpha_{3} , \alpha_{2} + \alpha_{3})$;
\item $a_{33} = 0$, $a_{11} a_{22} = a_{12} a_{21}$ and $a_{11} a_{33} = a_{13} a_{31}$; then $\beta_1 = \min (\alpha_{2} + \alpha_{3} , \alpha_{1})$, $\beta_3 = \max (\alpha_{1} + \alpha_{3} , \alpha_{2} + \alpha_{3})$, $\beta_2 = 2 \alpha_{1} + \alpha_{2} + 2 \alpha_{3}  - \beta_{1} - \beta_{3}$;
\item $a_{33} = 0$, $a_{11} a_{22} = a_{12} a_{21}$ and $a_{22} a_{33} = a_{23} a_{32}$; then $\beta_1 = \min (\alpha_{1} + \alpha_{3} , \alpha_{2})$, $\beta_3 = \max (\alpha_{1} + \alpha_{3} , \alpha_{2} + \alpha_{3})$, $\beta_2 = \alpha_{1} + 2 \alpha_{2} + 2 \alpha_{3}  - \beta_{1} - \beta_{3}$;
\item $a_{22} a_{33} \!=\! a_{23} a_{32}$, $a_{11} a_{33} \!=\! a_{13} a_{31}$ 
	and $a_{11} a_{22} \!=\! a_{12} a_{21}$; then $\beta_3 = \alpha_{j_2} + \alpha_{j_3}$, 
	$\beta_2 = \alpha_{j_1} + \alpha_{j_3}$, $\beta_1 = \alpha_{j_1} + \alpha_{j_2}$.
\end{enumerate}
In case \#1 we can additionally see that $a = a_{j_1,j_2} a_{j_2, j_1}$, $b = a_{j_1,j_3} a_{j_3, j_1}$ 
and $c = a_{j_2,j_3} a_{j_3, j_2}$; in case \#20 we have $a = a_{j_1, j_1}$, $b = a_{j_2, j_2}$ 
and $c = a_{j_3, j_3}$.
For the other cases, the computation of $a$, $b$ and $c$ from the given data in the differential
equation system~\eqref{eq:linearsystem} is also possible but requires very many case distinctions
that lead to general formulas which we believe to be far too unwieldy to write down here.
We have therefore decided to only write down the procedure for one specific example and believe
that the readers will be able to adapt the approach for the other cases when necessary:

\begin{exmp}
	\label{ex0a}
	Let us consider the case where $\alpha_1 = 0.4$, $\alpha_2 = 0.3$, $\alpha_3 = 0.5$
	and 
	\[
		A = \begin{pmatrix}
			-3 & 0 & 1.5 \\
			-0.5 & 0 & 0.5 \\
			6 & -1 & -3    
			\end{pmatrix} .
	\]
	Then we are in case \#15, and the characteristic function is
	\[
		Q(s) = s^{1.2} + 3 s^{0.8} + 3 s^{0.7} + 0.5 s^{0.4} + 0.75.
	\]
	Therefore, we can identify the parameters occurring in eq.~\eqref{eq:char-func-simple} as 
	$\beta_1 = \alpha_1 = 0.4$, $\beta_2 = \alpha_1 + \alpha_2 = 0.7$,
	$\beta_3 = \alpha_2 + \alpha_3 = 0.8$, $a = a_{11} = -3$, $b = a_{33} = -3$ and 
	$c = a_{22} a_{33} - a_{23} a_{32} = -0.5$. Moreover, $d = \det A = -0.75$. 
\end{exmp}

Our first positive main result then is the following sufficient condition for the asymptotic stability
of our differential equation system:
\begin{thm}
	\label{thm:stability1}
	If $Q(s)$ has the form \eqref{eq:char-func-simple} with 
	$d = \det A< 0$, $a , b , c \leq 0$, $\beta_1, \beta_2, \beta_3 \le 1$, 
	$\beta_4 - \beta_3 < 1$, $|\beta_1 + \beta_3 - \beta_4| \le 1$ and
	$|\beta_2 + \beta_3 - \beta_4| \le 1$
	then all roots of $Q(s)$ satisfy $ \Re(s) < 0$, 
	so that the system is asymptotically stable.
\end{thm}

\begin{pf}
	If $Q(s)$ has a root $s$ with $\Re(s) \geq 0$ then $| \arg (s) | \leq \frac{\pi}{2}$, and therefore
	$| \arg(s^{\beta_{i}}) | = \beta_{i} \cdot | \arg(s) | \leq \frac{\beta_{i} \pi}{2} \le \frac{\pi}{2}$
	for $i \in \{{1,2,3}\}$.
	In other words, we see that
	\begin{equation}
		\label{eq:re-gt-0}
		\Re(s^{\beta_{i}}) \geq 0 \mbox{ for } i \in \{1,2,3\}. 
	\end{equation}
	Moreover, since $Q(s) = 0$ 
	where $Q(s)$ is given by~\eqref{eq:char-func-simple}, we have
	\begin{equation*}
		s^{\beta_{3}} = \frac{b s^{\beta_{2}} + c s^{\beta_{1}} + d}{s^{\beta_{4} - \beta_{3}} - a}
	\end{equation*}
	where we have exploited the fact that, by our assumptions, $|\arg(s^{\beta_4 - \beta_3})| = 
	(\beta_4 - \beta_3) |\arg s| < \frac \pi 2$, so that $\Re (s^{\beta_4 - \beta_3}) > 0$ which---together
	with our assumption that $a \le 0$---implies that 
	the denominator on the right-hand side cannot be zero.
	Thus we obtain
	\begin{eqnarray*}
	\lefteqn{\Re(s^{\beta_{3}})} \\
	& = & \Re \left( \frac{b s^{\beta_{2}} + c s^{\beta_{1}} + d}{s^{\beta_{4} - \beta_{3}} - a} \right) \\
	& = &  \frac{ \Re \left( (b s^{\beta_{2}} + c s^{\beta_{1}} + d) (\bar{s}^{\beta_{4} - \beta_{3}} - a) \right)}{| s^{\beta_{4} - \beta_{3}} - a|^{2}}\\
	& = &  \frac 1 {| s^{\beta_{4} - \beta_{3}} - a|^{2}} \big[
		 \Re(b s^{\beta_{2}} + c s^{\beta_{1}} + d) \Re(\bar{s}^{\beta_{4} - \beta_{3}} - a) \\
	&& \qquad \qquad \qquad {} - \Im(b s^{\beta_{2}} + c s^{\beta_{1}} + d)\Im(\bar{s}^{\beta_{4} - \beta_{3}} - a)
		\big ]\\
	& = &  \frac 1 {| s^{\beta_{4} - \beta_{3}} - a|^{2}} \big[
				b (\Re(s^{\beta_2}) \Re(\bar s^{\beta_4 - \beta_3}) - \Im(s^{\beta_2}) \Im(\bar s^{\beta_4 - \beta_3}) ) \\
	&& \qquad \qquad \qquad {}
				- ab \Re(s^{\beta_2}) \! - \! ac \Re(s^{\beta_1}) + d \Re(\bar s^{\beta_4 - \beta_3})
				\! - \! ad  \\
	&& \qquad \qquad \qquad {} + \! c (\Re (s^{\beta_1}) \Re(\bar s^{\beta_4 - \beta_3})
					\! - \! \Im (s^{\beta_1}) \Im(\bar s^{\beta_4 - \beta_3}))
				\big] \\
	& = &  \frac 1 {| s^{\beta_{4} - \beta_{3}} - a|^{2}} \big[
				b \Re(s^{\beta_2} \bar s^{\beta_4 - \beta_3}) 
				- ab \Re(s^{\beta_2}) - ac \Re(s^{\beta_1}) \\
	& & \qquad \qquad \qquad {} + d \Re(\bar s^{\beta_4 - \beta_3})
				- ad +  c \Re (s^{\beta_1} \bar s^{\beta_4 - \beta_3})
				\big]				
	\end{eqnarray*}
	When we take into consideration the relation
	\[
		\Re (s^{\gamma_1} \bar s^{\gamma_2}) 
		= |s|^{\gamma_1 + \gamma_2} \cos ((\gamma_1 - \gamma_2) \arg s),
	\]
	the inequalities $a, b, c \leq 0$, $d < 0$, our assumptions on the $\beta_i$, and
	eq.~\eqref{eq:re-gt-0},
	we can see that the expression in the square brackets of the last member of the
	equation chain above is strictly negative. Since the factor outside of the square brackets is positive,
	we see that $\Re(s^{\beta_3}) < 0$ which contradicts our initial assumption above.
	So the equation $Q(s) = 0$ does not have any roots with $\Re(s) \geq 0$.
	\qed
\end{pf}

\begin{rem}
	Whether the conditions on the $\beta_i$ in this Theorem are fulfilled depends on which 
	of the 20 cases mentioned above we are in. In case 1, this is particularly simple to check because 
	then we have $\beta_4 = \beta_1 + \beta_2 + \beta_3 = \alpha_1 + \alpha_2 + \alpha_3$,
	so---since we always assume $\alpha_i \in (0,1]$---the assumptions are satisfied whenever
	$\alpha_{j_1} + \alpha_{j_2} < 1$.
\end{rem}

\begin{figure}[hbt]
\centering
	\includegraphics[width=0.75\columnwidth]{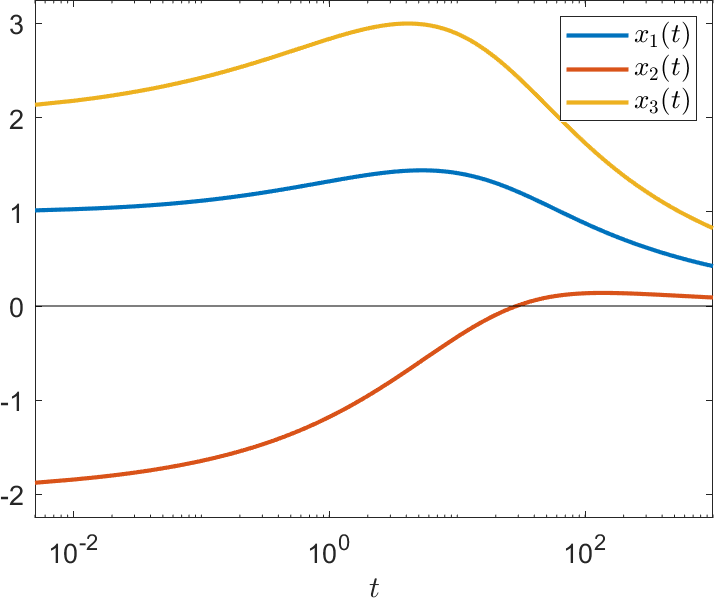}
	\caption{\label{fig:ex0}Solution to the system $D^\alpha x(t) = A x(t)$,
		$x(0) = (1, -2, 2)^{\mathrm T}$, with $A$ and $\alpha$ as in Example \ref{ex0a}. Note
		the logarithmic scale on the horizontal axis.}
\end{figure}

\begin{rem}
	\label{rem:ex0}
	It is easy to check that the problem discussed in Example \ref{ex0a} satisfies
	the conditions of Theorem \ref{thm:stability1}, so we can conclude that the 
	system $D^\alpha x(t) = A x(t)$ with $A$ and $\alpha$ as given there
	is asymptotically stable. We have plotted the components of the associated 
	numerical solution up to the value $t = 1000$ (computed with the implicit product integration rule of trapezoidal type 
	of \cite{Ga2018} with step size 1/200) 
	with initial values $x(0) = (1, -2, 2)^{\mathrm T}$ in Figure \ref{fig:ex0}.
	To display the main property of the solution components, viz.\ their
	convergence to zero for $t \to \infty$, as clearly as possible, the plots
	are shown in a coordinate system with a logarithmic scale on the horizontal axis.
\end{rem}

To find additional sets of conditions under which we can guarantee the zeros of the characteristic function
$Q$ to be in the open left half of $\mathbb C$, we need a few preparations. 
Remember that $Q$ has the form \eqref{eq:char-func-simple}.
For $\omega \in \mathbb R$, we write
\begin{align}
	\label{eq:def-h1}
	h_{1}(\omega)  =  \Re Q(\omega \mathrm i) 
	&  = \omega^{\beta_{4}} \cos\frac{\beta_{4} \pi}{2} 
		  - a \omega^{\beta_{3}} \cos\frac{\beta_{3} \pi}{2} \\
	& \phantom{==} {} - b \omega^{\beta_{2}} \cos\frac{\beta_{2} \pi}{2}  
		 - c \omega^{\beta_{1}} \cos\frac{\beta_{1} \pi}{2} -d \nonumber 
\intertext{and}	
	\label{eq:def-h2}
	h_{2}(\omega)  =  \Im Q(\omega \mathrm i) 
	& =    \omega^{\beta_{4}} \sin\frac{\beta_{4} \pi}{2} 
				 - a \omega^{\beta_{3}} \sin\frac{\beta_{3} \pi}{2} \\
	& \phantom{==}  {}   - b \omega^{\beta_{2}} \sin\frac{\beta_{2} \pi}{2} \nonumber 
		 - c \omega^{\beta_{1}} \sin\frac{\beta_{1} \pi}{2}.
\end{align}

First, we look at zeros of $Q$ on the imaginary axis. 

\begin{lem}
	\label{lem:beta4}
	Let $d <0$.
	If $Q$ has a zero on the imaginary axis then there exists some $\omega > 0$ such that 
	$h_1(\omega) = h_2(\omega) = 0$.
\end{lem}

\begin{pf}
	Since $Q(\bar s) = \overline{Q(s)}$,
	it follows that $Q(\omega  \mathrm i) = 0$ with some $\omega > 0$ if and only if 
	$Q(- \omega \mathrm i) = 0$. Moreover, from $d < 0$ we see that $Q(0) \ne 0$. 
	Thus, the statement follows directly from the definitions 
	\eqref{eq:def-h1} and \eqref{eq:def-h2} of $h_1$ and $h_2$.
%
	\qed
\end{pf}

Next, to investigate the existence of zeros of $Q$ in the open right half of $\mathbb C$ 
we recall an auxiliary result already developed in \citet[pp.~1332--1333]{DiethelmThaiTuan}.
In view of Theorem \ref{thm:unstable1}, it is clear that such zeros exist if $d > 0$ and that 
a zero at the origin exists if $d = 0$. It thus only remains to consider the case $d < 0$.
Here we define the sets
\begin{align*}
	\gamma_1 &:=\left \{s= \omega \mathrm i: \varepsilon \leq \omega \leq R  \right \}, \\
	\gamma_2 &:=\left \{ s=R e^{\mathrm i \varphi}: -\frac{\pi}{2}\leq \varphi  \leq \frac{\pi}{2} \right\}, \\
	\gamma_3 &:=\left \{ s=\varepsilon e^{\mathrm i \varphi}: -\frac{\pi}{2}\leq \varphi  \leq \frac{\pi}{2} \right\}, \\
	\gamma_4 &:=\left \{s= \omega \mathrm i: -R \leq \omega \leq -\varepsilon  \right \}
\end{align*}
where $0 < \varepsilon < R$ are chosen such that all zeros of $Q$ in the open right half of $\mathbb C$ are
located inside of the open region $C_{\varepsilon,R} \subset \mathbb C$ whose boundary is formed by these four curves. 
This is possible if we choose $\varepsilon$ sufficiently small (since $Q(0) > 0$, the continuity of $Q$ 
implies that $Q$ has no zeros $s$ with $|s| \le \varepsilon$) and $R$ sufficiently large (because 
$\lim_{|s| \to \infty} |Q(s)| = \infty$ uniformly in $\arg s$).
As in \citet[pp.~1332--1333]{DiethelmThaiTuan}, we see that the number of zeros of $Q$ in the 
open right half of $\mathbb C$ is equal to the number of zeros in $C_{\varepsilon,R}$ which, since $Q$ is analytic in 
$C_{\varepsilon,R}$, is equal (by the argument principle) to the number $Z$
of times that the curve $Q(\cup_{k=1}^4 \gamma_k)$ encircles the origin in the positive direction.
If $Z > 0$ then this curve must encircle the origin at least once, and so the curve must cross
the positive and negative real and imaginary axes. Thus, in particular there is some 
$s^{*} \in \cup_{k=1}^4 \gamma_{k}$ such that $Q(s^{*}) \in (-\infty , 0)$. 

When $\beta_4 < 2$, using the argument above, it is no loss of generality to assume that $\varepsilon$ and $R$ are such that  $s^{*} \notin \gamma_2 \cup \gamma_{3}$, so $s^{*} \in \gamma_{1} \cup \gamma_{4} \subset \{\omega \mathrm i : \omega \in \mathbb R\}$. Since $Q(\bar{s}) = \overline{Q(s)}$ 
there must exist some $\omega \in (\varepsilon, R) \subset \mathbb{R}$ such that $Q(\omega \mathrm i) \in (-\infty, 0)$.
We have thus shown:

\begin{lem}
	\label{lem:qzeros}
	If $\beta_4 < 2$ and $Q$ has zeros in the open right half of $\mathbb C$ then there exists an $\omega > 0$ 
	with $h_1(\omega) < 0$ and $h_2(\omega) = 0$.
\end{lem}

Hence, in view of Lemmas \ref{lem:beta4} and \ref{lem:qzeros}, 
we must find conditions under which $h_{1}(\omega) > 0$ holds for all $\omega > 0$ with $h_{2}(\omega) = 0$.

To simplify the notation for our following developments, we now assume that $\beta_4 \ne 2$ and introduce the quantities
\begin{equation}
	\label{eq:rho}
	\rho_i = \frac{\sin \frac{(\beta_4 - \beta_i) \pi}2}{\sin \frac{\beta_4 \pi}2}
	\quad \text{ and } \quad 
	\tilde \rho_i = \frac{\sin \frac{\beta_i \pi}{2}}{\sin \frac{\beta_4 \pi}2}
\end{equation}
for $i = 1,2,3$, so that whenever $h_2(\omega) = 0$ we have
\begin{eqnarray}
	h_1(\omega) \nonumber
	&=& 	\frac {h_1(\omega) \sin \frac{\beta_4 \pi} 2 
			- h_2(\omega) \cos \frac{\beta_4 \pi} 2} {\sin \frac{\beta_4 \pi} 2  } \\
	&=& - a \omega^{\beta_{3}} \rho_3 - b \omega^{\beta_2} \rho_2 - c \omega^{\beta_1} \rho_1 - d.
	\label{s}
\end{eqnarray}

%
%
%

\begin{lem}
	\label{lem:tt}
	Let $0 \! < \! \beta_1 \! < \! \beta_2 \! < \! \beta_3 \! < \! \beta_4 \! < \! 2$. If $a,b,c \le 0$ and $d < 0$
	then all zeros of $Q$ are in the open left half of $\mathbb C$.
\end{lem}

\begin{pf}
	Since $\rho_k > 0$ for $k = 1,2,3$, our assumptions on $a$, $b$, $c$ and $d$ imply in view of 
	\eqref{s} that $h_1(\omega) > 0$ for all $\omega > 0$ with $h_2(\omega) = 0$.
	The conclusion then follows by Lemmas \ref{lem:beta4} and \ref{lem:qzeros}.
	\qed
\end{pf}

\begin{lem}   
	\label{lem1}
	If, in eq.~\eqref{eq:char-func-simple}, $a = b = 0$, $c > 0$, $\beta_4 < 2$ and
	\begin{equation}
		\label{a}
		d < - c \rho_{1} (c \tilde{\rho_{1}})^{\beta_{1}/(\beta_{4}-\beta_{1})},
	\end{equation}
	then all zeros of $Q$ are in the open left half of $\mathbb C$.
\end{lem}

\begin{pf}    
	In view of our assumptions, we see that $h_{2}(\omega_{0}) = 0$ if and only if 
	$\omega_{0} = ( c \tilde{\rho_{1}})^{1 / (\beta_{4} - \beta_{1})}$. 
	From (\ref{s}), we have 
	$ h_{1}(\omega_{0}) = - c (c \tilde{\rho_{1}})^{\beta_{1}/(\beta_{4}-\beta_{1})} \rho_{1} - d $,
	so \eqref{a} implies $ h_{1}(\omega_{0}) > 0 $. 
	Thus, by Lemmas \ref{lem:beta4} and \ref{lem:qzeros}, 
	$Q$ has no zeros on the imaginary axis or in the open right half of $\mathbb C$.
	\qed
\end{pf}

In a completely analog manner, we can also show:

\begin{lem}
	\label{lem2}
	Assume that, in eq.~\eqref{eq:char-func-simple}, $\beta_4 < 2$ and
	\begin{itemize}
	\item[(a)] either $b = c= 0$, $a > 0$ and 
		$d < - a \rho_{3} (a\tilde{\rho_{3}})^{\beta_{3}/(\beta_{4}-\beta_{3})}$
	\item[(b)] or $a = c = 0$, $b > 0$ and 
		$d < - b \rho_{2} (b\tilde{\rho_{2}})^{\beta_{2}/(\beta_{4}-\beta_{2})}$.
	\end{itemize}
	Then, all zeros of $Q$ lie in the open left half of $\mathbb C$.
\end{lem}

\begin{lem}
	\label{lem4}
	Let $\beta_1 < \beta_2 < \beta_{3}$ and $\beta_4 < 2$. Assume
	\begin{enumerate}
	\item[(i)] either $a = 0$, $b > 0$, $c > 0$,
		 $1 < b \tilde \rho_{2} + c \tilde \rho_{1}$ and
			\[
			        d \le - b ( b \tilde \rho_2 + c \tilde{\rho_{1}})^\frac{\beta_{2}}{\beta_{4} - \beta_{2}} \rho_{2} 
					  - c ( b \tilde \rho_2 + c \tilde{\rho_{1}})^\frac{\beta_{1}}{\beta_{4} - \beta_{2}} \rho_{1};
			\]
	\item[(ii)] or  $a > 0$, $b = 0$, $c > 0$,
		$1 < a \tilde \rho_{3} + c \tilde \rho_{1}$ and
			\[
			        d \le - a ( a \tilde \rho_3 + c \tilde{\rho_{1}})^\frac{\beta_{3}}{\beta_{4} - \beta_{3}} \rho_{3} 
					  - c ( a \tilde \rho_3 + c \tilde{\rho_{1}})^\frac{\beta_{1}}{\beta_{4} - \beta_{3}} \rho_{1};
			\]
	\item[(iii)] or  $a > 0$, $b > 0$, $c = 0$,
		$1 < a \tilde \rho_{3} + b \tilde \rho_{2}$ and
			\[
			        d \le - a ( a \tilde \rho_3 + b \tilde{\rho_{2}})^\frac{\beta_{3}}{\beta_{4} - \beta_{3}} \rho_{3} 
					  - b ( a \tilde \rho_3 + b \tilde{\rho_{2}})^\frac{\beta_{2}}{\beta_{4} - \beta_{3}} \rho_{2}.
			\]
	\end{enumerate}
	Then all zeros of $Q$ are in the open left half of $\mathbb C$.
\end{lem}

\begin{pf}
	We have 
	\begin{align*}
		h_{2}'(\omega) & =  \omega^{\beta_{1} - 1} g_{2}(\omega) \\
	\intertext{where}
		g_{2}(\omega)  
		& =   \beta_{4} \omega^{\beta_{4} - \beta_{1}} \sin\frac{\beta_{4} \pi}{2} - a \beta_{3} \omega^{\beta_{3} - \beta_{1}} \sin\frac{\beta_{3} \pi}{2} \\
		& \phantom{==}
			{ } - b \beta_{2} \omega^{\beta_{2} - \beta_{1}} \sin\frac{\beta_{2} \pi}{2} 
			- c \beta_{1} \sin\frac{\beta_{1} \pi}{2}.
	\end{align*}

	Now for part (i), since $a = 0$, $ b > 0$ and $c > 0$, we have 
	\begin{eqnarray*}
		g_{2}(\omega) 
		& = & \omega^{\beta_{2} - \beta_{1}} 
			\left( \beta_{4} \omega^{\beta_{4} - \beta_{2}} \sin\frac{\beta_{4} \pi}{2} -  b \beta_{2} \sin\frac{\beta_{2} \pi}{2} \right) \\ 
		&& {}  - c \beta_{1} \sin\frac{\beta_{1} \pi}{2} \\
		&=& c_{1} \omega^{d_{1}} + c_{2} \omega^{d_{2}} + c_{3}
	\end{eqnarray*}
	where 
	$c_{1} = \beta_{4} \sin (\beta_{4} \pi/2) > 0$, $c_{2} = -b \beta_{2} \sin (\beta_{2} \pi / 2) < 0$, 
	$c_{3} = -c \beta_{1} \sin (\beta_{1} \pi / 2) < 0$, 
	$d_{1} = \beta_{4} - \beta_{1}$ and $d_{2} = \beta_{2} - \beta_{1}$.
	It is not difficult to check that $d_{1} > d_2 > 0$. Thus, since
	\[
		{g}_{2}'(\omega)  
		\! = \! c_1 d_1 \omega^{d_1 - 1} + c_2 d_2 \omega^{d_2 - 1}
		\! = \! \omega^{d_2 - 1} (c_1 d_1 \omega^{d_1 - d_2} + c_2 d_2 ),
	\]
	we can see that $g_2'(\omega) > 0$ for 
	\[ 
		\omega > \omega_1 := \left( - \frac{c_2 d_2}{c_1 d_1} \right)^{1/(d_1-d_2)}
	\]
	and $g_2'(\omega) < 0$ for $\omega \in (0, \omega_1)$. Due to the facts that 
	$g_{2}(0) = c_3 < 0$ if $\beta_2 > \beta_1$, $g_2(0) = c_2 + c_3 < 0$ if $\beta_2 = \beta_1$,
	and $\lim_{\omega\rightarrow +\infty} g_{2}(\omega) = + \infty$, $g_{2}$ has exactly one root in $(0, \infty)$.
	Therefore, $h_2'$ also has exactly one root in $(0, \infty)$. Furthermore, $\lim_{\omega \to 0} h_2(\omega) = 0$ and
	$\lim_{\omega \to \infty} h_2(\omega) = \infty$. Thus, $h_2$ also has exactly one root $\omega_0 \in (0, \infty)$.

	Since $b \tilde \rho_{2} + c \tilde \rho_{1} > 1$, we have $h_{2}(1) < 0$, so $\omega_0 > 1$. 
	Moreover, due to $\beta_{3} > \beta_{2} > \beta_{1} > 0$, we have 
	\begin{eqnarray*}
		h_{2}(\omega) 
		& > & \omega^{\beta_{4}} \sin\frac{\beta_{4} \pi}{2}
			 - \left(b \sin\frac{\beta_{2} \pi}{2} + c  \sin\frac{\beta_{1} \pi}{2} \right) \omega^{\beta_{2}} \\
		& = & \omega^{\beta_{2}} \sin\frac{\beta_{4} \pi}{2} 
				\left ( \omega^{\beta_4 - \beta_2} - \left(b \tilde\rho_2 + c \tilde \rho_1 \right) \right)
	\end{eqnarray*}
	for every $\omega > 1$, and thus $h_{2} ((b \tilde{\rho_{2}} + c \tilde{\rho_{1}})^{1/(\beta_{4} - \beta_{2})}) > 0$. 
	This implies that $1 < \omega_{0} < (b \tilde{\rho_{2}} + c \tilde{\rho_{1}})^{1/(\beta_{4} - \beta_{2})}$. Hence, 
	in view of eq.~\eqref{s} and our assumption on the value of $d$,
	\begin{eqnarray*}
		h_1(\omega_0) 
		& = & - b \omega_0^{\beta_2} \rho_2 - c \omega_0^{\beta_1} \rho_1 - d \\
		& > & - b  (b \tilde{\rho_{2}} + c \tilde{\rho_{1}})^{\beta_2/(\beta_{4} - \beta_{2})} \rho_2 \\
		&&		{} - c  (b \tilde{\rho_{2}} + c \tilde{\rho_{1}})^{\beta_1/(\beta_{4} - \beta_{2})} \rho_1 - d \\
		& \ge & 0
	\end{eqnarray*}
	which completes the proof in case (i).
	
	For cases (ii) and (iii), we proceed in the same way.
	\qed
\end{pf}

\begin{pf*}{\bf Proof of Theorem \ref{thm:unstable1}(b).}
	As in \eqref{eq:def-h2}, set $h_2(\omega) = \Im Q(\omega \mathrm i)$.
	Then, since $\psi_k \in (0,2)$, $\sin \frac{\psi_k \pi}2 \ge 0$ 
	for all $k$. Moreover, $\beta_4 = \alpha_1 + \alpha_2 + \alpha_3 \le 3$, and thus
	our assumption that $\beta_4 \ge 2$ leads to
	$\sin \frac{\beta_4 \pi}2 \le 0$. From $\max c_k > \min c_k \ge 0$ we then derive 
	$h_2(\omega) = \Im Q(\omega \mathrm i) < 0$ for all $\omega > 0$.

	Now we return to the considerations leading up to the proof of Lemma \ref{lem:qzeros}.
	Since $Q(0) = -\det A > 0$, we may assume $\varepsilon$ to be sufficiently small so that 
	$\Re Q(s) >0$ for all $s \in \gamma_3$, Therefore, the curve $Q(\gamma_3)$ is entirely
	contained in the open right half of $\mathbb C$. 
	Furthermore, since $h_2(\omega) < 0$ for $\omega > 0$, the curve
	$Q(\gamma_1)$ is completely contained in the open bottom half of $\mathbb C$, 
	and by the symmetry relation $Q(\bar s) = \overline{Q(s)}$, $Q(\gamma_4)$  lies in the
	open top half of $\mathbb C$. 

	If the value $R$ used in the definition of $\gamma_2$
	is sufficiently large (which we may assume without loss of generality) then 
	$Q(s)$ essentially behaves as $s^{\beta_4}$ on the semicircle $\gamma_2$. Thus,
	as $\varphi$ evolves from $-\pi/2$ to $\pi/2$ when we proceed through the curve $\gamma_2$,
	$\arg Q(s)$ increases from $-\beta_4 \pi/2 < -\pi$ to $\beta_4 \pi/2 > \pi$. Combining this
	with the other observations above, we see that $Q(s)$ encirlces the origin twice in counterclockwise 
	direction, and therefore $Q$ has two zeros inside the region $C_{\varepsilon,R}$ and hence
	in the open right half of $\mathbb C$, so Theorem \ref{thm:deng} implies the instability.
	\qed
\end{pf*}

\subsection{Inhomogeneous linear equations}
\label{subs:inhom-lin}

For inhomogeneous linear systems with certain types of forcing functions $f_k$,
one can prove that the stability properties can be fully investigated solely on the basis of the stability 
of the associated homogeneous system. The fundamental result in this context reads as follows, see
\citet[Theorem~5.3]{DiethelmHashemishahrakiThaiTuan}.

\begin{thm}
	\label{thm:inhom-lin}
	Consider the initial value problem \eqref{eq:linearsystem} with $\alpha_k \in (0,1]$ for all $k$ and
	an arbitrarily chosen initial condition $x_k(0) = x_k^0 \in \mathbb R$ ($k = 1, 2, 3$). 
	Denote the forcing function vector by
	$f(t) = (f_1(t), f_2(t), f_3(t))^{\mathrm T}$,
	set $\mu = \min \{\alpha_1, \alpha_2, \alpha_3\}$ and assume that 
	all zeros of the associated characteristic function $Q$ (as defined in Theorem \ref{thm:deng})
	are in the open left half of the complex plane. Then we have:
	\begin{itemize}
	\item [(i)] If $f$ is bounded then the solution of the initial value problem is also bounded.
	\item [(ii)] If $\lim_{t \to \infty} f(t) = 0$ then the solution of the initial value problem 
		converges to $0$ as $t \to \infty$.
	\item [(iii)] If $\|f(t) \| = O(t^{-\eta})$ as $t \to \infty$ with some $\eta > 0$
		then the solution $x$ of the initial value problem
		satisfies $\|x(t)\| = O(t^{-\nu}) $ as $t \to \infty$ where $\nu = \min\{\mu, \eta\}$.  
	\end{itemize}
\end{thm}

It is thus essentially unnecessary to deal with the inhomogeneous case separately; all required stability results
can normally be obtained by investigating the corresponding properties of the associated homogeneous system.

\subsection{Autonomous nonlinear problems}
\label{subs:nonlin}

The approach introduced for the two-dimensional case in \citet[Theorem 5]{DiethelmThaiTuan}
allows to also apply our findings to certain autonomous nonlinear differential equation systems. However, 
while Theorem \ref{thm:inhom-lin} deals with inhomogeneous linear systems
for arbitrary initial values, we now have to restrict ourselves 
to initial values that are sufficiently close to the equilibrium solution.
The main result reads as follows \cite[Theorem~5.5]{DiethelmHashemishahrakiThaiTuan}.

\begin{thm}
	\label{thm:nonlin}
	Consider the initial value problem
	\begin{equation}
		\label{eq:ivp-nonlin}
		\left.
		\begin{array}{rl}
			D^\alpha x(t) & =  A x(t) + f(x(t)), \\
			x(0) & =  x^0,
		\end{array}
		\right \}
	\end{equation}
	where $\alpha = (\alpha_1, \alpha_2, \alpha_3) \in (0,1]^3$, $A$ is a constant
	$(3 \times 3)$ matrix, and $x^0 \in \mathbb R^3$.
	Let the function $f: \Omega \to \mathbb R^3$ with $0 \in \Omega \subset \mathbb R^3$
	satisfy $f(0) = 0$ and fulfil a local Lipschitz condition in a neighbourhood of the origin
	such that
	\[
		\lim_{r\to 0} \sup_{0 \le \|x\|, \|y\| < r,\;x\ne y}\frac{\|f(x)-f(y)\|}{\|x-y\|} = 0.
	\]
	If all zeros of the characteristic function $Q$ associated to the matrix $A$
	as in Theorem \ref{thm:deng} are in the open left half of $\mathbb C$, 
	there exists $\delta > 0$ such that, whenever $\|x^0\| < \delta$,
	the unique global solution $x$ of the problem~\eqref{eq:ivp-nonlin} satisfies
	$\| x(t) \| = O(t^{-\nu})$ as $t \to \infty$
	with $\nu = \min \{\alpha_1, \alpha_2, \alpha_3\}$.
\end{thm}

\section{A numerical example}
\label{sec:examples}

We now take a look at a specfic example to demonstrate the applicability of our results.

\begin{exmp}
	\label{ex1}
	Consider the system
	\begin{equation}
		\label{eq:ex1}
		\left.
		\begin{array}{ll}
	    		D^{0.4} x_{1}(t) & =  x_{2}(t) - x_{3}(t) + f_{1} (t), \\
	    		D^{0.3} x_{2}(t) & = 0.2 x_{1}(t) + f_{2}(t), \\
	    		D^{0.5} x_{3}(t) & = 0.5 x_{2}(t) + f_{3}(t),
	  	\end{array}
		\right \}
	\end{equation}
	with
	\begin{equation*}
		f_{k}(t) 
		=	\begin{cases}
    				1           & \text{ if } 0 \leq t < 1, \\  
	    			t^{-2 k}  & \text{ if } t \geq 1
  			\end{cases}
	\end{equation*}
	for $k = 1, 2, 3$.
	These forcing functions satisfy the conditions of 
	each of the three parts of Theorem \ref{thm:inhom-lin} with $\eta = 2$ in part (iii).
	The characteristic function of the homogeneous system is 
	$Q(s) = s^{1.2} - 0.2 s ^{0.5} + 0.1$, and thus
	all zeros of $Q$ lie in the open left half of the complex plane by Lemma~\ref{lem1}.
	Therefore, Theorem \ref{thm:inhom-lin} allows us
	to deduce the asymptotic stability of the inhomogeneous system.
\end{exmp}

\begin{figure}[ht]
	    \centering
	    \includegraphics[width=0.75\columnwidth]{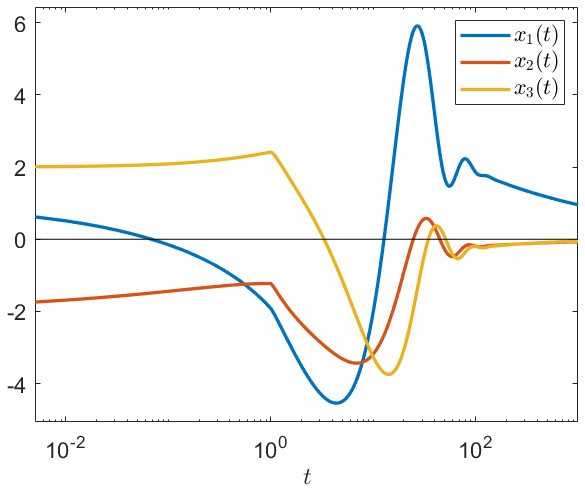}
	    \caption{Solution to differential equation~\eqref{eq:ex1} from Example~\ref{ex1}
	    		with initial values $(1, -2, 2)^{\mathrm T}$. The horizontal axis is plotted in a logarithmic scale.}
	    \label{fig:example1a}
	\end{figure}

	The asymptotic stability of the solution is evident from Figure \ref{fig:example1a}
	that shows a plot of the components of the solution on the interval  $[0, 1000]$. 
	
\begin{figure}[ht]
	    \centering
	    \includegraphics[width=0.75\columnwidth]{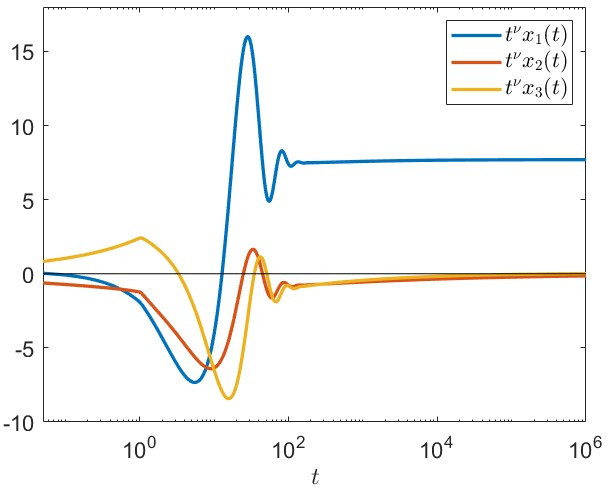}
	    \caption{The functions $t^{\nu} x_{j}(t)$ with $\nu = 0.3$ where the $x_j$ are the solutions to the
	    	differential equation~\eqref{eq:ex1} from Example~\ref{ex1}
	  	with initial values $(1, -2, 2)^{\mathrm T}$. The horizontal axis has a logarithmic scale. 
	  	The asymptotic behavior predicted by
	  	Theorem \ref{thm:inhom-lin}(iii) is visible.}
	    \label{fig:example1b}
	\end{figure}

	Figure \ref{fig:example1b} exhibits the decay properties of the solution components
	as stipulated by Theorem \ref{thm:inhom-lin}(iii). It turns out that we need to look at a much longer
	interval than before to really see the boundedness of $t^\mu x(t)$, especially for the first component.
	Note that, in this example, $\eta = 2$ (in the notation of Theorem \ref{thm:inhom-lin}), so 
	$\nu = \min \{ \alpha_1, \alpha_2, \alpha_3, \eta \} = 0.3$.
	
\section{Conclusion}

We have presented some necessary and some sufficient conditions for the asymptotic stability of a three-dimensional incommensurate homogeneous linear differential system with Caputo derivatives, and we have demonstrated how these results can be applied to obtain corresponding conditions for inhomogeneous linear systems and for certain types of autonomous nonlinear systems. Numerical examples are provided to illustrate the main theoretical results.

A number of questions remain open for the moment, e.g.\ the case $\beta_4 \ge 2$ 
in Lemmas \ref{lem1}, \ref{lem2} and \ref{lem4} and the case when the coefficients of the matrix do not 
satisfy the conditions of any of the cases (1)--(20) in Subsection \ref{subs:lin-hom}. 
We intend to address these issues in our future work.

%
%
%



\bibliography{root}             
                                                   







\end{document}